\documentclass[a4paper,draft]{amsproc}
\usepackage{amssymb}
\usepackage[hyphens]{url} 

\theoremstyle{plain}
 \newtheorem{thm}{Theorem}[section]
 
 \newtheorem{lem}{Lemma}[section]
 
\theoremstyle{definition}
 
 \newtheorem{dfn}{Definition}[section]
\theoremstyle{remark}
 
 \numberwithin{equation}{section}

\renewcommand{\leq}{\leqslant}
\renewcommand{\geq}{\geqslant}

\title[New variations of power increasing sequences]{New variations of power increasing sequences}

\subjclass[2010]{26D15; 42A24; 40F05; 40G99}

\keywords{Summability factors, absolute matrix summability, Fourier series, infinite series, H\"{o}lder inequality, Minkowski inequality}

\author[YILDIZ]{\bfseries \c{S}ebnem Y{\i}ld{\i}z} 

\address{Department of Mathematics \\
Ahi Evran University\\ K{\i}r\c{s}ehir, Turkey}
\email{sebnemyildiz@ahievran.edu.tr; sebnem.yildiz82@gmail.com}

\begin{document}

{\begin{flushleft}\baselineskip9pt\scriptsize
\end{flushleft}}
\vspace{18mm} \setcounter{page}{1} \thispagestyle{empty}

\begin{abstract}
The aim of this paper is to generalize a main theorem concerning weighted mean summability to absolute matrix summability which plays a vital role in  summability theory and applications to the other sciences by using quasi-$f$-power sequences.
\end{abstract}

\maketitle

\section{Introduction}  

\begin{dfn}\cite{Bar} A positive sequence $(b_{n})$ is said to be an almost increasing sequence if there exists a positive increasing sequence $(c_{n})$ and two positive constants $M$ and $N$ such that $Mc_{n}\leq b_{n}\leq Nc_{n}$.
\end{dfn}
\begin{dfn}\cite{Su3}	
	A positive sequence $X=(X_{n})$ is said to be quasi-$f$-power increasing sequence if there exists a constant $K=K(X,f)\geq 1$ such that $Kf_{n}X_{n}\geq f_{m}X_{m}$ for all $n\geq m\geq 1$, where $f=\left\lbrace f_{n}(\sigma,\beta)\right\rbrace =\left\lbrace n^{\sigma}(logn)^{\beta},\beta\geq 0, 0<\sigma<1 \right\rbrace$.
\end{dfn}
\begin{dfn}
 The sequence $(\lambda_{n})$ is said to be of bounded variation, denoted by $(\lambda_{n})\in \mathcal{BV}$, if $\sum\limits_{n=1}^{\infty}|\Delta\lambda_{n}|<\infty.$
 If we take $\beta=0$, then we have a quasi-$\sigma$-power increasing sequence. Every almost increasing sequence is a quasi-$\sigma$-power increasing sequence for any non-negative $\sigma$, but the converse is not true for $\sigma>0$ (see \cite{Le}). For any sequence $(\lambda_{n})$ we write that $\Delta^{2}\lambda_{n}=\Delta\lambda_{n}-\Delta\lambda_{n+1}$ and $\Delta\lambda_{n}=\lambda_{n}-\lambda_{n+1}$.
\end{dfn}
Let $\sum a_{n}$ be a given infinite series with the partial sums $(s_{n})$. 
By $u_{n}^{\alpha}$ and $t_{n}^{\alpha}$ we denote the nth Ces\`{a}ro means of order $\alpha$, with $\alpha> -1$, of the sequence $(s_{n})$ and $(na_{n})$, respectively, that is (see \cite{EC})
\begin{align}
u_{n}^{\alpha}=\frac{1}{A_{n}^{\alpha}}\sum_{v=0}^{n} A_{n-v}^{\alpha-1}s_{v}\quad	\textnormal{and} \quad t_{n}^{\alpha}=	\frac{1}{A_{n}^{\alpha}}\sum_{v=0}^{n} A_{n-v}^{\alpha-1}va_{v}, 						
\end{align}
\noindent where
\begin{align}
A_{n}^{\alpha}=\frac{(\alpha +1)(\alpha +2)...(\alpha+n)}{n!}=O(n^{\alpha}),~~ \quad  A_{-n}^{\alpha}=0 \quad \textnormal{for} \quad n>0.
\end{align} 
\begin{dfn}\cite{TMF},\cite{EK}   
	The series $\sum a_{n}$ is said to be summable $|C,\alpha |_{k}$, $k\geq 1$, if
	\begin{align}
	\sum_{n=1}^{\infty}n^{k-1}|u_{n}^{\alpha}-u_{n-1}^{\alpha}|^{k}=\sum_{n=1}^{\infty}\frac{1}{n}|t_{n}^{\alpha}|^{k}<\infty.
	\end{align}
\end{dfn}
 If we take $\alpha=1$, then $|C, \alpha|_{k}$ summability reduces to $|C,1|_{k}$ summability.\\
Let $(p_{n})$ be a sequence of positive real numbers such that
\begin{align}
P_{n}=\sum_{v=0}^{n}p_{v}\rightarrow \infty \quad{as}
\quad{n}\rightarrow \infty,\quad (P_{-i}=p_{-i}=0,\quad  i\geq 1).
\end{align}
The sequence-to-sequence transformation
\begin{align}
t_{n}=\frac{1}{P_{n}}\sum_{v=0}^{n}p_{v}s_{v}
\end{align}
defines the sequence $(t_n)$ of the Riesz mean or simply the
$(\bar{N},p_n)$ mean of the sequence $(s_{n})$ generated by the
sequence of coefficients $(p_n)$ (see \cite{GH}).
\begin{dfn}\cite{Bor1}
	The series $\sum{a_n}$
	is said to be summable $\left| \bar{N},p_{n}\right| _k$, $k \geq 1,$ if 
	\begin{align}
	\sum_{n=1}^{\infty}\left( \frac{P_{n}}{p_{n}}\right) ^{k-1}\mid t_{n}-t_{n-1}\mid^k
	< \infty.
	\end{align}
\end{dfn}
In the special case when $p_{n}=1$ for all values of $n$ (resp.~$k=1$), $\left| \bar{N},p_{n}\right| _k$  summability is the same as $\left|C,1\right|_{k} $ (resp.~$\mid{\bar{N},p_n}\mid$) summability.
\section {The Known Results}
The following theorems are known dealing with the $\left| \bar{N},p_{n}\right| _k$ summability factors of infinite series.
 \begin{thm}{\cite{Maz}} Let $(X_{n})$ be a almost increasing sequence. If the sequences $(X_{n})$, $(\lambda_{n})$, and $(p_{n})$ satisfy the conditions 
	\begin{align}\label{eq:2}
	\lambda_{m}X_{m}&=O(1)  \quad as \quad {m}\rightarrow \infty,\\
	\sum_{n=1}^{m}nX_{n}|\Delta^{2}\lambda_{n}|&=O(1) \quad as \quad {m}\rightarrow \infty,\\
	\sum_{n=1}^{m}\frac{P_{n}}{n}&=O(P_{m})\\
	\sum_{n=1}^{m}\frac{p_{n}}{P_{n}}|t_{n}|^{k}&=O(X_{m}) \quad  as \quad {m}\rightarrow \infty,\\
	\sum_{n=1}^{m}\frac{|t_{n}|^{k}}{n}&=O(X_{m})\quad  as \quad {m}\rightarrow \infty,
	\end{align}
	then the series $\sum a_{n}\lambda_{n}$ is summable $\left| \bar{N},p_n\right| _k$,  $k\geq 1$.
 \end{thm} 
\begin{thm}{\cite{Bor6}} Let $(X_{n})$ be a quasi-$\sigma$-power increasing sequence. If the sequences $(X_{n})$, $(\lambda_{n})$ and $(p_{n})$ satisfy the conditions (2.1)-(2.3), and 
	\begin{align}
	\sum_{n=1}^{m}\frac{p_{n}}{P_{n}}\frac{|t_{n}|^{k}}{X_{n}^{k-1}}&=O(X_{m}) \quad  as \quad {m}\rightarrow \infty,\\
	\sum_{n=1}^{m}\frac{|t_{n}|^{k}}{nX_{n}^{k-1}}&=O(X_{m})\quad  as \quad {m}\rightarrow \infty,
	\end{align}
	then the series $\sum a_{n}\lambda_{n}$ is summable $\left| \bar{N},p_{n}\right| _k$, $k\geq 1$.
\end{thm}
Later on, Bor has proved the following theorem by taking quasi-f-power increasing sequence instead of a quasi-$\sigma$-power increasing sequence.
\begin{thm}\textnormal{\cite{Bor7}} Let $(X_{n})$ be a quasi-$f$-power increasing sequence. If the sequences $(X_{n})$, $(\lambda_{n})$ and $(p_{n})$ satisfy all the conditions of Theorem 2.2,
	then the series $\sum a_{n}\lambda_{n}$ is summable $\left| \bar{N},p_{n}\right| _k$,  $k\geq 1$.
\end{thm}
Let $A=(a_{nv})$ be a normal matrix,
i.e., a lower triangular matrix of nonzero diagonal entries. Then $A$ defines the sequence-to-sequence transformation, mapping the
sequence $s=(s_{n})$ to $As=\left(A_{n}(s)\right)$, where
\begin{align}\label{eq:6}
A_{n}(s)=\sum_{v=0}^{n}a_{nv}s_{v}, \quad n=0,1,...
\end{align}

\begin{dfn}\cite{WT}
	The series $\sum a_{n}$ is said to be summable $\left| A,p_{n} \right|_{k}$, $k\geq 1$, if 
	\begin{align}\label{eq:7}
	\sum_{n=1}^{\infty}\left(\frac{P_{n}}{p_{n}} \right) ^{k-1}\left|\bar{\Delta}A_{n}(s)\right|^{k}< \infty,
	\end{align}
	where
	\begin{align}\label{eq:9}
	\bar{\Delta}A_{n}(s)=A_{n}(s)-A_{n-1}(s).
	\end{align}
\end{dfn}
If we take $p_{n}=1$ for all values of $n$, then we have $\left| A\right|_{k}$ summability (see \cite{NT}). And also if we take  $a_{nv}=\frac{p_{v}}{P_{n}}$, then we have $\left| \bar{N},p_n\right| _k$ summability. Furthermore, if we take $a_{nv}=\frac{p_{v}}{P_{n}}$ and $p_{n}=1$ for all values of $n$, then $\left| A,p_{n} \right|_{k}$ summability reduces to $\left| C,1\right|_{k}$ summability (see \cite{TMF}).
\section {The Main Results}
The Fourier series play an important role in many areas of applied mathematics and mechanics. Recently some papers have been done concerning absolute matrix summability of infinite series and Fourier series (see \cite{Bor2}-\cite{Bor5}, \cite{oz1}-\cite{Sa}, \cite{Y3}-\cite{Y5}). 
The aim of this paper is to generalize Theorem 2.3 for $|A,p_{n}|_{k}$ summability method for these series by taking quasi-f-power increasing sequence instead of a quasi-$\sigma$-power increasing sequence.\\
Given a normal matrix $A=(a_{nv})$, we associate two lower
semimatrices $\bar{A}=(\bar{a}_{nv})$ and $\hat{A}=(\hat{a}_{nv})$
as follows:
\begin{align}\label{eq:11}
\bar{a}_{nv}=\sum_{i=v}^{n}a_{ni},\quad n,v=0,1,...
\end{align}
and
\begin{align}\label{eq:12}
\hat{a}_{00}=\bar{a}_{00}=a_{00},\quad
\hat{a}_{nv}=\bar{a}_{nv}-\bar{a}_{n-1,v},\quad n=1,2,...
\end{align}
It may be noted that $\bar{A}$ and $\hat{A}$ are the well-known
matrices of series-to-sequence and series-to-series
transformations, respectively. Then, we have
\begin{align}\label{eq:13}
A_{n}(s)&=\sum_{v=0}^{n}a_{nv}s_{v}=
\sum_{v=0}^{n}\bar{a}_{nv}a_{v}
\end{align} 
and
\begin{align}
\bar{\Delta}A_{n}(s)&=\sum_{v=0}^{n}\hat{a}_{nv}a_{v}.
\end{align}
Using this notation we have the following theorem.
 \begin{thm}\label{thm:4} Let $(X_{n})$ be a quasi-$f$-power increasing sequence. Let $k\geq 1$ and $A=(a_{nv})$ be a positive normal matrix such that
	\begin{align}\label{eq:15}
	\overline{a}_{no}&=1,\     n=0,1,...,\\
	a_{n-1,v}&\geq a_{nv},\ \textnormal{for}~~   n\geq v+1,\\
	a_{nn}&=O\left(\frac{p_{n}}{P_{n}} \right) \\
	\sum_{v=1}^{n-1}\frac{1}{v}\hat{a}_{n,v+1}&=O(a_{nn}).
	\end{align}
	If the sequences $(X_{n})$, $(\lambda_{n})$ and $(p_{n})$ satisfy all the conditions of Theorem 2.3, then the series $ \sum a_{n}\lambda_{n}$ is summable $\left|A,p_{n}\right|_{k}$, $k\geq 1$.
\end{thm}

It may be remarked that if we take $A=(\bar{N},p_{n})$, the conditions $(3.5)$-$(3.7)$ are satisfied automatically and the condition $(3.8)$ is satisfied by the condition $(2.3)$.
We need the following lemmas for the proof of our theorem.
\begin{lem}\label{lem:1} \textnormal{\cite{Bor2}} Under the conditions of  Theorem 2.1 we have that
	\begin{align}\label{eq:20}
	nX_{n}|\Delta\lambda_{n}|&=O(1)\quad as \quad  n\to \infty,\\
	\sum_{n=1}^{\infty}X_{n}|\Delta\lambda_{n}|&<\infty.
	\end{align}
\end{lem} 
\section*{Proof of Theorem~\ref{thm:4}}
	Let $X_{n}$ be a be a quasi-$f$-power increasing sequence and  $(I_{n})$ denotes the A-transform of the series $\sum_{n=1}^{\infty} a_{n} \lambda_{n}$. Then, we have
	\begin{align*}
	\bar{\Delta}I_{n}& =  \sum_{v=1}^{n}\hat{a}_{nv}a_{v}\lambda_{v}.
	\end{align*}
	Applying Abel's transformation to this sum, we have that
	\begin{align*}
	\bar{\Delta}I_{n}& =  \sum_{v=1}^{n}\hat{a}_{nv}a_{v}\lambda_{v}\frac{v}{v}=\sum_{v=1}^{n-1} \Delta (\frac{\hat{a}_{nv}\lambda_{v}}{v}) \sum_{r=1}^{v} r a_{r} + \frac{\hat{a}_{nn}\lambda_{n}}{n} \sum_{r=1}^{n} r a_{r}\\	
	& =  \sum_{v=1}^{n-1} \Delta(\frac{\hat{a}_{nv}\lambda_{v}}{v})(v+1)t_{v}+\hat{a}_{nn}\lambda_{n}\frac{n+1}{n}t_{n}\nonumber\\	
	& =  \sum_{v=1}^{n-1} \bar{\Delta}a_{nv}\lambda_{v}t_{v}\frac{v+1}{v}+\sum_{v=1}^{n-1}\hat{a}_{n,v+1}\Delta \lambda_{v}t_{v}\frac{v+1}{v}+ \sum_{v=1}^{n-1}\hat{a}_{n,v+1}\lambda_{v+1}\frac{t_{v}}{v}+a_{nn}\lambda_{n}t_{n}\frac{n+1}{n}\nonumber\\
	& =  I_{n,1}+I_{n,2}+I_{n,3}+I_{n,4}.
	\end{align*}
	To complete the proof of Theorem~\ref{thm:4}, by Minkowski's inequality, it is sufficient to show that
	\begin{align}\label{eq:20}
	\sum_{n=1}^{\infty}\left(\frac{P_{n}}{p_{n}} \right) 
	^{k-1}\mid
	I_{n,r}\mid^{k}< \infty, \quad \textnormal{for} \quad{r=1,2,3,4.}
	\end{align}
	First, by applying H{\"o}lder's inequality with indices \textit{k} and
	$k'$, where $\textit{k}>1$
	and $\frac{1}{k}+\frac{1}{k'}=1$, we have that
	\begin{align*}
	&	\sum_{n=2}^{m+1}\left(\frac{P_{n}}{p_{n}} \right) 
	^{k-1}\mid I_{n,1}\mid^{k}
	\leq \sum_{n=2}^{m+1}\left(\frac{P_{n}}{p_{n}} \right)
	^{k-1}\left\lbrace \sum_{v=1}^{n-1}|\frac{v+1}{v}|\left|\bar{\Delta}a_{nv}\right|
	|\lambda_{v}||t_{v}|
	\right\rbrace ^{k}\\
	&= O(1) \sum_{n=2}^{m+1}\left(\frac{P_{n}}{p_{n}} \right) 
	^{k-1}\sum_{v=1}^{n-1}\left|\bar{\Delta}a_{nv}\right||\lambda_{v}|^{k}|t_{v}|^{k} \times \left\lbrace     \sum_{v=1}^{n-1}\left|\bar{\Delta}a_{nv}\right| \right\rbrace ^{k-1},
	\end{align*}
	using
	\begin{align*}
	\Delta\hat{a}_{nv}&=\hat{a}_{nv}-\hat{a}_{n,v+1}=\bar{a}_{nv}-\bar{a}_{n-1,v}-\bar{a}_{n,v+1}+\bar{a}_{n-1,v+1}=a_{nv}-a_{n-1,v},
	\end{align*}
	and from $(3.5)$ and $(3.6)$ we have
	\begin{align*} 
	\sum_{v=1}^{n-1}|\bar{\Delta}a_{nv}|&=\sum_{v=1}^{n-1}|a_{nv}-a_{n-1,v}|=\sum_{v=1}^{n-1}(a_{n-1,v}-a_{nv})\\
	&=\sum_{v=0}^{n-1}a_{n-1,v}-a_{n-1,0}-\sum_{v=0}^{n}a_{nv}+a_{n0}+a_{nn}\\
	&=1-a_{n-1,0}-1+a_{n0}+a_{nn}\leq a_{nn},
	\end{align*}
	and using $\sum_{n=v+1}^{m+1}|\bar{\Delta}a_{nv}|\leq a_{vv}$ 
	we have,   	
	\begin{align*}		
	&\sum_{n=2}^{m+1}\left(\frac{P_{n}}{p_{n}} \right) 
	^{k-1}\mid I_{n,1}\mid^{k}= O(1)
	\sum_{n=2}^{m+1}\left(\frac{P_{n}}{p_{n}} \right) 
	^{k-1}a_{nn}^{k-1}\left\lbrace  \sum_{v=1}^{n-1}|\bar{\Delta}a_{nv}||\lambda_{v}|^{k}|t_{v}|^{k}\right\rbrace\\
	&	= O(1)
	\sum_{v=1}^{m}|\lambda_{v}|^{k-1}|\lambda_{v}||t_{v}|^{k}\sum_{n=v+1}^{m+1}|\bar{\Delta}a_{nv}| \nonumber \\ 							  
	&= O(1)
	\sum_{v=1}^{m}\frac{1}{X_{v}^{k-1}}|\lambda_{v}||t_{v}|^{k}a_{vv}= O(1)
	\sum_{v=1}^{m-1}\Delta|\lambda_{v}|\sum_{r=1}^{v}a_{rr}\frac{|t_{r}|^{k}}{X_{r}^{k-1}}+O(1)|\lambda_{m}|
	\sum_{v=1}^{m}a_{vv}\frac{|t_{v}|^{k}}{X_{v}^{k-1}}\nonumber\\	
	&=O(1)
	\sum_{v=1}^{m-1}|\Delta\lambda_{v}|X_{v}+O(1)|\lambda_{m}|X_{m}=O(1)\quad \textnormal{as}\quad m\rightarrow\infty,
	\end{align*}
	by virtue of the hypotheses of  Theorem~\ref{thm:4} and Lemma 3.1. Also, we have that
	\begin{align*}
	&\sum_{n=2}^{m+1}\left(\frac{P_{n}}{p_{n}} \right) ^{k-1}\mid
	I_{n,2}\mid^{k} \leq
	\sum_{n=2}^{m+1}\left(\frac{P_{n}}{p_{n}} \right) 
	^{k-1}\left\lbrace \sum_{v=1}^{n-1}|\frac{v+1}{v}||\hat{a}_{n,v+1}||\Delta\lambda_{v}||t_{v}|\right\rbrace^{k}\nonumber \\
	& = O(1)\sum_{n=2}^{m+1}\left(\frac{P_{n}}{p_{n}} \right) 
	^{k-1}\left\lbrace \sum_{v=1}^{n-1}\hat{a}_{n,v+1}|\Delta\lambda_{v}||t_{v}|\frac{X_{v}}{X_{v}}\right\rbrace ^{k}\nonumber\\
	& = O(1)\sum_{n=2}^{m+1}\left(\frac{P_{n}}{p_{n}} \right) 
	^{k-1}\left\lbrace \sum_{v=1}^{n-1}\hat{a}_{n,v+1}|\Delta\lambda_{v}|X_{v}\frac{1}{X_{v}^{k}}|t_{v}|^{k}\right\rbrace \times \left\lbrace\sum_{v=1}^{n-1}\hat{a}_{n,v+1}|\Delta\lambda_{v}|X_{v}\right\rbrace^{k-1}\\ 
	& = O(1)\sum_{n=2}^{m+1}\left(\frac{P_{n}}{p_{n}} \right) 
	^{k-1}a_{nn}^{k-1}\left\lbrace\sum_{v=1}^{n-1}\hat{a}_{n,v+1}|\Delta\lambda_{v}|X_{v}\frac{1}{X_{v}^{k}} |t_{v}|^{k}  \right\rbrace \times \left\lbrace\sum_{v=1}^{m-1}|\Delta\lambda_{v}|X_{v} \right\rbrace ^{k-1}\\
	& = O(1)\sum_{v=1}^{m}v|\Delta\lambda_{v}|\frac{1}{X_{v}^{k-1}}\frac{1}{v}|t_{v}|^{k}\sum_{n=v+1}^{m+1}\hat{a}_{n,v+1}= O(1)\sum_{v=1}^{m}v|\Delta\lambda_{v}|\frac{1}{vX_{v}^{k-1}}|t_{v}|^{k}\\ 
	& = O(1)\sum_{v=1}^{m-1}\Delta(v|\Delta\lambda_{v}|)\sum_{r=1}^{v}\frac{|t_{r}|^{k}}{rX_{r}^{k-1}}+O(1)m|\Delta\lambda_{m}|\sum_{r=1}^{m}\frac{|t_{r}|^{k}}{rX_{r}^{k-1}}\\
	& = 	O(1)\sum_{v=1}^{m-1}|\Delta(v|\Delta\lambda_{v}|)|X_{v}+
	O(1)m|\Delta\lambda_{m}|X_{m}\nonumber \\      
	& = 
	O(1)\sum_{v=1}^{m-1}vX_{v}|\Delta^{2}\lambda_{v}|+O(1)\sum_{v=1}^{m-1}X_{v}|\Delta\lambda_{v}|+O(1)m|\Delta\lambda_{m}|X_{m}\nonumber \\
	& = 
	O(1)\quad \textnormal{as}\quad m\rightarrow\infty,
	\end{align*}
	by virtue of the hypotheses of Theorem~\ref{thm:4} and Lemma 3.1.  
	Furthermore, as in $I_{n,1}$, we have
	\begin{align*}
	&	\sum_{n=2}^{m+1}\left(\frac{P_{n}}{p_{n}} \right) 
	^{k-1}\mid I_{n,3}\mid^{k}\leq
	\sum_{n=2}^{m+1}\left(\frac{P_{n}}{p_{n}} \right) 
	^{k-1}\left\lbrace \sum_{v=1}^{n-1}|\hat{a}_{n,v+1}| |\lambda_{v+1}|\frac{|t_{v}|}{v}\right\rbrace ^{k} \nonumber\\
	& = 
	O(1)\sum_{n=2}^{m+1}\left(\frac{P_{n}}{p_{n}} \right) 
	^{k-1}\left\lbrace\sum_{v=1}^{n-1}
	|\hat{a}_{n,v+1}| |\lambda_{v+1}|^{k}\frac{|t_{v}|^{k}}{v}\right\rbrace \times\left\lbrace \sum_{v=1}^{n-1}\frac{1}{v}\hat{a}_{n,v+1}\right\rbrace ^{k-1}\nonumber \\	
	& = 
	O(1)\sum_{n=2}^{m+1}\left(\frac{P_{n}}{p_{n}} \right) 
	^{k-1}a_{nn}^{k-1}\sum_{v=1}^{n-1}|\lambda_{v+1}||\lambda_{v+1}|^{k-1}\frac{|t_{v}|^{k}}{v}\hat{a}_{n,v+1}= 
	O(1)\sum_{v=1}^{m}\frac{|t_{v}|^{k}}{v}\frac{1}{X_{v}^{k-1}}|\lambda_{v+1}|
	\sum_{n=v+1}^{m+1}\hat{a}_{n,v+1}\nonumber \\
	& = 
	O(1)\sum_{v=1}^{m}\frac{|t_{v}|^{k}}{v}\frac{1}{X_{v}^{k-1}}|\lambda_{v+1}|
	\sum_{n=v+1}^{m+1}\hat{a}_{n,v+1}= 
	O(1)\sum_{v=1}^{m}\frac{1}{X_{v}^{k-1}}|\lambda_{v+1}|\frac{|t_{v}|^{k}}{v}
	\nonumber \\
	& = 
	O(1)\quad \textnormal{as} \quad m\rightarrow\infty,
	\end{align*}
	by virtue of the hypotheses of  Theorem 3.1 and Lemma 3.1. Again, as in $I_{n,1}$, we have that
	\begin{align*}
	\sum_{n=1}^{m}\left(\frac{P_{n}}{p_{n}} \right) 
	^{k-1}|I_{n,4}|^{k} &=O(1) \sum_{n=1}^{m}\left(\frac{P_{n}}{p_{n}} \right) 
	^{k-1}a_{nn}^{k}|\lambda_{n}|^{k}|t_{n}|^{k}= O(1)\sum_{n=1}^{m}\left(\frac{P_{n}}{p_{n}} \right) 
	^{k-1}a_{nn}^{k-1}a_{nn} |\lambda_{n}|^{k-1}|\lambda_{n}||t_{n}|^{k}\\
	&=O(1)\sum_{n=1}^{m}a_{nn}\frac{1}{X_{n}^{k-1}}|\lambda_{n}||t_{n}|^{k}= O(1)\quad \textnormal{as} \quad m\rightarrow\infty,
	\end{align*}
	by virtue of hypotheses of the Theorem~\ref{thm:4} and Lemma 3.1.
	This completes the proof of Theorem 3.1.
\section{An application of absolute matrix summability to Fourier series}
Let $f$ be a periodic function with period $2\pi$ and integrable $(L)$ over $(-\pi,\pi)$. Without any loss of generality  the constant term in the Fourier series of $f$ can be taken to be zero, so that
\begin{align}
f(t)\sim \sum_{n=1}^{\infty}(a_{n}cosnt+b_{n}sinnt)=\sum_{n=1}^{\infty}C_{n}(t).
\end{align}
where
\begin{align}
a_{0}&=\frac{1}{\pi}\int_{-\pi}^{\pi}f(t)dt,\quad\nonumber
a_{n}=\frac{1}{\pi}\int_{-\pi}^{\pi}f(t)cos(nt)dt, \nonumber\quad
b_{n}=\frac{1}{\pi}\int_{-\pi}^{\pi}f(t)sin(nt)dt.\nonumber
\end{align}
We write
\begin{align}
\phi(t)=\frac{1}{2}\left\lbrace f(x+t)+f(x-t)\right\rbrace,\\
\phi_{\alpha}(t)= \frac{\alpha}{t^{\alpha}} \int_{0}^{t} (t-u)^{\alpha-1}\phi(u)\, du, \quad (\alpha> 0).
\end{align}
It is well known that if $\phi(t)\in \mathcal{BV}(0,\pi)$, then $t_{n}(x)=O(1)$, where  $t_{n}(x)$ is the $(C,1)$ mean of the sequence $(nC_{n}(x))$ (see \cite{KKC}).\\ Using this fact, Bor has obtained the following main result dealing with the trigonometric Fourier series.

\begin{thm}\cite{Bor6}
	Let $(X_{n})$ be a quasi-$\sigma$-power increasing sequence. If $\phi_{1}(t)\in \mathcal{BV}(0,\pi)$, and the sequences $(p_{n})$, $(\lambda_{n})$, and $(X_{n})$ satisfy the conditions of Theorem 2.3, then the series $\sum C_{n}(x)\lambda_{n}$ is summable $|\bar{N},p_{n}|_{k}$, $k\geq 1$.	
\end{thm}
\begin{thm}\cite{Bor7}
	Let $(X_{n})$ be a quasi-$f$-power increasing sequence. If $\phi_{1}(t)\in \mathcal{BV}(0,\pi)$, and the sequences $(p_{n})$, $(\lambda_{n})$, and $(X_{n})$ satisfy the conditions of Theorem 2.3, then the series $\sum C_{n}(x)\lambda_{n}$ is summable $|\bar{N},p_{n}|_{k}$, $k\geq 1$.	
\end{thm}
We now apply the above theorems to the weighted mean in which $A=(a_{nv})$ is defined as $a_{nv}=\frac{p_{v}}{P_{n}}$ when $0\leq v\leq n$, where $P_{n}=p_{0}+p_{1}+...+p_{n}.$
Therefore, it is well known that \begin{align*}
\bar{a}_{nv}=\frac{P_{n}-P_{v-1}}{P_{n}} \quad \text{and}\quad \hat{a}_{n,v+1}=\frac{p_{n}P_{v}}{P_{n}P_{n-1}}. 
\end{align*}
So, one can easily verify that the conditions of Theorem 3.1 reduce to those of Theorem 2.3 and also we can obtain new results dealing with absolute matrix summability of Fourier series in the following manner. 
\begin{thm}
	Let $A$ be a positive normal matrix satisfying the conditions of Theorem 3.1. Let $(X_{n})$ be a quasi-$\sigma$-power increasing. If $\phi_{1}(t)\in \mathcal{BV}(0,\pi)$, and the sequences $(p_{n})$, $(\lambda_{n})$, and $(X_{n})$ satisfy the conditions of Theorem 3.1, then the series $\sum C_{n}(x)\lambda_{n}$ is summable $|A,p_{n}|_{k}$, $k\geq 1$.	
\end{thm}
\begin{thm}
	Let $A$ be a positive normal matrix satisfying the conditions of Theorem 3.1. Let $(X_{n})$ be a quasi-$f$-power increasing sequence. If $\phi_{1}(t)\in \mathcal{BV}(0,\pi)$, and the sequences $(p_{n})$, $(\lambda_{n})$, and $(X_{n})$ satisfy the conditions of Theorem 3.1, then the series $\sum C_{n}(x)\lambda_{n}$ is summable $|A,p_{n}|_{k}$, $k\geq 1$.	
\end{thm}
\section{\small{APPLICATIONS}}
We may now ask whether there are some examples other than weighted mean methods of matrices $A$ that satisfy the hypotheses of Theorem 3.1. For example, apply Theorem 3.1 to the Ces\`{a}ro method of order $\alpha$ with $0<\alpha\leq 1$ in which $A$ is given by $a_{nv}=A_{n-v}^{\alpha-1}/A_{n}^{\alpha}$, and by applying 
Theorem 3.1, Theorem 4.3 and Theorem 4.4 to weighted mean so, the following results can be easily verified.\\
\noindent1. If we take $a_{nv}=\frac{p_{v}}{P_{n}}$ in Theorem 3.1, Theorem 4.3 and Theorem 4.4, then we have Theorem 2.3, Theorem 4.1 and Theorem 4.2.\\
2. If we take $\beta=0$ and $a_{nv}=\frac{p_{v}}{P_{n}}$ in Theorem 3.1 and Theorem 4.4, then we have Theorem 2.2 and Theorem 4.1.\\
3.  If we take $p_{n}=1$ for all values of $n$ in Theorem 3.1, Theorem 4.3 and Theorem 4.4, then we have a new  result dealing with $\left| A\right|_{k}$ summability.\\
4. If we take $a_{nv}=\frac{p_{v}}{P_{n}}$ and $p_{n}=1$ for all values of $n$ in Theorem 3.1, Theorem 4.3 and Theorem 4.4, then we have a new result concerning $\left| C,1\right|_{k}$ summability.

\end{document}